\documentclass[12pt]{article}
\usepackage{amsmath, amssymb, amsfonts, amsthm}
\usepackage{color}
\newcommand{\R}{\mathbb R}
\newcommand{\C}{\mathbb C}

\newcommand{\be}{\begin{equation}}
\newcommand{\ee}{\end{equation}}

\def\0{\mathbf 0}

\newtheorem{thm}{Theorem}[section]

\newtheorem{corollary}[thm]{Corollary}

\theoremstyle{remark}
\newtheorem{rem}[thm]{Remark}
\def \bui#1#2{\mathrel{\mathop{\kern 0pt#1}\limits^{#2}}}
\numberwithin{equation}{section}

\newcommand{\eq}{\begin{equation}}
\newcommand{\eeq}{\end{equation}}
\usepackage[lmargin=2.5 cm,rmargin=2.5 cm,tmargin=3.5cm,bmargin=2.5cm,paper=a4paper]{geometry}
\begin{document}
\title{\bf Eigenvalue bounds of the Robin Laplacian with magnetic field}
\author{Georges Habib\footnote{Lebanese University, Faculty of Sciences II, Department of Mathematics, P.O. Box 90656 Fanar-Matn, Lebanon,
E-mail: \texttt{ghabib@ul.edu.lb}},\, Ayman Kachmar \footnote{Lebanese University, Faculty of Sciences V, Department of Mathematics, Nabatieh, Lebanon, E-mail: \texttt{ayman.kashmar@gmail.com}}}

\date{}
\maketitle

\begin{abstract}
\noindent On a compact Riemannian manifold $M$ with boundary, we give an estimate for the eigenvalues $(\lambda_k(\tau,\alpha))_k$ of the magnetic Laplacian with the Robin boundary conditions. Here, $\tau$ is a positive number that defines the Robin condition and $\alpha$ is a real differential 1-form on $M$ that represents the magnetic field. We express these estimates in terms of the mean curvature of the boundary, the parameter $\tau$  and a lower bound of the Ricci curvature of $M$ (see Theorem \ref{estimate1} and Corollary \ref{corestimate}). The main technique is to use the Bochner formula established in \cite{ELMP} for the magnetic Laplacian and to integrate it over $M$ (see Theorem \ref{bochnermagnetic1}). In the last part, we compare the eigenvalues $\lambda_k(\tau,\alpha)$ with the first eigenvalue $\lambda_1(\tau)=\lambda_1(\tau,0)$ (i.e. without magnetic field) and the Neumann eigenvalues $\lambda_k(0,\alpha)$ (see Theorem \ref{thm:comp}) using the min-max principle.
\end{abstract}

%\noindent{\it Key words}: Magnetic Laplacian, Robin condition, magnetic field, eigenvalue.\\

%\noindent{\it Mathematics Subject Classification}: ......

\section{Introduction and Results}
Let $(M,g)$ be a Riemannian manifold of dimension $n$ and let $\alpha$ be a smooth real differential $1$-form on $M.$ Given two vector fields $X,Y$ in  the complexified tangent bundle $TM\otimes \C,$ the {\it magnetic covariant derivative} is defined as $\nabla_Y^\alpha X=\nabla^M_YX+i\alpha(Y)X,$ where $\nabla^M$ denotes the Levi-Civita connection on $M.$  It is shown in \cite[Lemma 3.2]{ELMP} that $\nabla^\alpha$ satisfies the Leibniz rule and the  compatibility property with respect to the Riemannian metric $g$, and is also used to define the {\it magnetic Hessian} by ${\rm Hess}^\alpha f(X,Y)=\langle \nabla_X^\alpha d^\alpha f,Y\rangle.$ Here and in all the paper, the product $\langle\cdot,\cdot\rangle$ will denote the Hermitian inner product extended from the metric $g$ to the tangent bundle $TM\otimes \mathbb{\C}$ or to the cotangent bundle $T^*M\otimes \mathbb{\C}.$ We will also use the natural one-to-one isomorphism between $T^{*}M\otimes \C$ and $TM\otimes \C$ by $w(X)=\langle X,\overline{w^{\#}}\rangle$ for any $X\in TM\otimes \C$ and $w\in T^{*}M\otimes \C.$
 
\noindent
 Given any complex-valued function $f$ on $M,$ the {\it magnetic Laplacian} is defined as being the trace of the magnetic Hessian
$$\Delta^\alpha f:=-{\rm trace}({\rm Hess}^\alpha f)=-{\rm div}^\alpha(d^\alpha f)^{\#},$$
where $d^\alpha f:=d^M f+if\alpha$ and ${\rm div}^\alpha$ is the magnetic divergence given for any vector field $X\in TM\otimes \C$ by ${\rm div}^\alpha X:={\rm div}^M X+i\langle X,\alpha^{\#}\rangle.$ 

\noindent
 The study of the spectrum of the magnetic Laplacian has interested many researchers \cite{CS, Er, FLM, Sh, Shu1, Shu2} during the last years. For example, the authors in \cite{ELMP} gave an estimate {\it \`a la Lichnerowicz} for the first eigenvalue in terms of a lower bound of the Ricci curvature (assumed to be positive) and the infinity norm of the magnetic field $d^M\alpha$. In particular, they deduce a spectral gap between the first eigenvalue (which is not necessarily zero) and the second one. The main technique used in the paper is a Bochner type formula for the magnetic Laplacian $\Delta^\alpha,$ which they integrate it over the manifold $M$ and they control all the integral terms involving $d^M\alpha$. Indeed, they prove

\begin{thm} \label{bochnermagnetic} \cite[Thm. 4.1]{ELMP} Let $(M,g)$ be a complete Riemannian manifold of dimension $n.$ Then for all $f\in C^\infty(M,\C)$, we have
\begin{eqnarray} \label{eq:bochner} 
-\frac{1}{2}\Delta^M(|d^\alpha f|^2)&=&|{\rm Hess}^\alpha f|^2-\Re\langle d^\alpha f,d^\alpha(\Delta^\alpha f)\rangle+{\rm Ric}^M(d^\alpha f,d^\alpha f)\nonumber\\&+&i(d^M\alpha(d^\alpha f,\overline{d^\alpha f})-d^M\alpha(\overline{d^\alpha f},d^\alpha f))\nonumber\\&+&\frac{i}{2}(\langle \bar{f}d^\alpha f,\delta^M d^M\alpha\rangle-\langle f\overline{d^\alpha f},\delta^M d^M\alpha\rangle),\nonumber\\
\end{eqnarray}
where $\delta^M$ denotes the formal adjoint of $d^M$ on $(M,g).$ 
\end{thm} 

\noindent In this paper, we are interested in estimating the eigenvalues of the magnetic Laplacian with the Robin boundary condition. That is, we assume on a given compact manifold $M$ with boundary $N$ there exists a complex-valued function $f$ on $M$ satisfying the equation $\Delta^\alpha f=\lambda f$ on $M$ and the boundary condition $(d^\alpha f)(\nu)=\tau f$ for some positive real number $\tau.$ Here $\nu$ denotes the inward unit normal vector field of $N,$ which will be identified with its dual one form. It a standard fact that the spectrum of such boundary problem is purely discrete and consists of a sequence of eigenvalues  $(\lambda_k(\tau,\alpha))_k$ arranged in increasing order counting multiplicities. In order to get the estimates for the eigenvalues, we shall first integrate the Bochner formula in Theorem \ref{bochnermagnetic} as in \cite{ELMP} by taking into account the boundary terms. First, we get

\begin{thm} \label{bochnermagnetic1} Let $(M^n,g)$ be a compact Riemannian manifold with boundary $N$ and let $\alpha$ be a differential real $1$-form on $M.$ Then, we have
\begin{align}\label{bochnermagnetic2}
\int_M|{\rm Hess}^\alpha f+\frac{1}{n}(\Delta^\alpha f) g|^2 dv_g=&
\frac{n-1}{n}\int_M|\Delta^\alpha f|^2dv_g-\int_M {\rm Ric}^M(d^\alpha f,d^\alpha f) dv_g\nonumber\\&+\int_M \Im m\left((d^M\alpha)(d^\alpha f,\overline{d^\alpha f})\right)dv_g+\int_M|f|^2|d^M\alpha|^2dv_g\nonumber\\&-(n-1)\int_N H|\langle d^\alpha f,\nu\rangle|^2dv_g
-2\int_N\Re (\langle\nu,d^\alpha f\rangle\Delta_N^\alpha f)dv_g\nonumber\\&-\int_N\langle II(d^\alpha_Nf),d_N^\alpha f\rangle dv_g.
\end{align}
for all complex valued function $f\in C^\infty(M,\mathbb{C}).$ 
\end{thm}

\noindent Here $II$ denotes the second fundamental form of the boundary and $H$ is the mean curvature. Also $\Delta_N^\alpha$ is a Laplacian defined on functions on $N$ which is associated to some exterior derivative $d^\alpha_N$ (see Section \ref{sec:proofbochner} for the definition). 

\noindent The formula \eqref{bochnermagnetic2} can be useful for different applications in spectral theory. One of these applications is to use Theorem \ref{bochnermagnetic1} for a particular solution of the magnetic Robin boundary problem. Therefore, we get the universal bound on the eigenvalues of the magnetic Robin Laplacian under some assumptions on the magnetic field $d^M\alpha$, the Ricci curvature ${\rm Ric}^M$ and the second fundamental form $II.$ Indeed,
 
\begin{thm} \label{estimate1} Let $(M^n,g)$ be a compact Riemannian manifold with boundary $\partial M=N$ and let $\alpha$ be a differential 1-form on $M$ and $\tau>0.$ Assume that ${\rm Ric}^M\geq k\,(k >0)$ and that $II+\tau\geq 0.$ If $\alpha$ satisfies 
\begin{equation} \label{conditionalpha}
k-(n-1)\tau H_{\min}\leq ||d^M\alpha||_\infty\leq \left(1+2\sqrt{\frac{n-1}{n}}\right)^{-1}k,
\end{equation}
then any eigenvalue $\lambda(\tau,\alpha)$ of the Laplacian $\Delta^\alpha$ satisfies
$$\lambda(\tau,\alpha)\leq a_-(k,||d^M\alpha||_\infty,n) \quad\text{or}\quad \lambda(\tau,\alpha)\geq a_+(k,||d^M\alpha||_\infty,n),$$
where 
$$a_\pm(k,||d^M\alpha||_\infty,n)=n\frac{(k-||d^M\alpha||_\infty)\pm\sqrt{(k-||d^M\alpha||_\infty)^2-4(\frac{n-1}{n})||d^M\alpha||_\infty^2}}{2(n-1)},$$
and $H_{\min}:=\min_M H$.
\end{thm}

\begin{rem}\label{rem:estimate1}~
\begin{itemize}
\item The assumption in \eqref{conditionalpha} on the  mean curvature is valid when $H_{\min}> 0$, since $\left(1+2\sqrt{\frac{n-1}{n}}\right)^{-1}k<k$.
Also, when $\tau$ is very large, \eqref{conditionalpha} becomes an upper bound on $||d^M\alpha||_\infty$, which is a growth condition on the magnetic field with respect to the Ricci curvature. 

\item It follows from Inequality \eqref{conditionalpha} that $(k-||d^M\alpha||_\infty)^2-4(\frac{n-1}{n})||d^M\alpha||_\infty^2>0$ and
$a_-(k,\|d^M\alpha\|_\infty,n)> 0$. This is more transparent in the proof of Theorem \ref{estimate1}.
\end{itemize}
\end{rem} 

\noindent As a direct consequence of Theorem \ref{estimate1} and a standard continuity argument as in \cite{ELMP}, one gets
  
\begin{corollary} \label{corestimate} Let $(M^n,g)$ be a compact Riemannian manifold with boundary $\partial M=N$ and let $\alpha$ be a differential 1-form on $M$ and $\tau>0.$ Assume that ${\rm Ric}^M\geq k\,(k >0)$ and that $II+\tau\geq 0.$ If $k\leq (n-1)\tau H_{\min}$ and $\alpha$ satisfies 
\begin{equation*} 
||d^M\alpha||_\infty\leq \left(1+2\sqrt{\frac{n-1}{n}}\right)^{-1}k,
\end{equation*}
then any eigenvalue $\lambda(\tau,\alpha)$ of the Laplacian $\Delta^\alpha$ satisfies
$$\lambda(\tau,\alpha)\geq a_+(k,||d^M\alpha||_\infty,n),$$
where 
$$a_+(k,||d^M\alpha||_\infty,n)=n\frac{(k-||d^M\alpha||_\infty)+\sqrt{(k-||d^M\alpha||_\infty)^2-4(\frac{n-1}{n})||d^M\alpha||_\infty^2}}{2(n-1)}.$$
\end{corollary}

\noindent {\it Proof of Corollary  \ref{corestimate}:} It is enough to prove the lower bound on the first eigenvalue $\lambda_1(\tau,\alpha)$. We apply Theorem \ref{estimate1} to the $1$-form $\alpha'=\varepsilon \alpha,$ for $\varepsilon \in ]0,1[.$ The inequality \eqref{conditionalpha} is clearly satisfied for $\alpha'.$ Hence  $\lambda_1(\tau,\varepsilon\alpha)$ is either less than $a_-(k,\varepsilon||d^M\alpha||_\infty,n)$ or bigger than $a_+(k,\varepsilon||d^M\alpha||_\infty,n).$ Note that $\lambda_1(\tau,\varepsilon\alpha)$ and $a_-(k,\varepsilon||d^M\alpha||_\infty,n)$ depend continuously on $\varepsilon$. Since $\lambda_1(\tau,0)>0$ and  $a_-(k,\varepsilon||d^M\alpha||_\infty,n)\mathop{\longrightarrow}\limits_{\varepsilon\rightarrow 0}0$, we get that the inequality $\lambda_1(\tau,\varepsilon\alpha)\geq a_+(k,\varepsilon||d^M\alpha||_\infty,n)$ is true in a neighborhood of $\varepsilon=0$. Define $\varepsilon_*=\sup\{\varepsilon\in\,(0,1)~|~ \lambda_1(\tau,\varepsilon\alpha)\geq a_+(k,\varepsilon||d^M\alpha||_\infty,n)\}.$ If $\varepsilon_*<1$, then we get $\lambda_1(\tau,\varepsilon_*\alpha)\geq  a_+(k,\varepsilon_*||d^M\alpha||_\infty,n)$ and $\lim\limits_{\delta\to 0_+}\lambda_1(\tau,(\varepsilon_*+\delta)\alpha)\leq  a_-(k,\varepsilon_*||d^M\alpha||_\infty,n)$, which violates the continuity of $\lambda_1(\tau,\varepsilon\alpha)$ with respect to $\varepsilon$. Therefore, $\varepsilon_*=1$.
\hfill$\square$\\

\noindent As a direct application of Corollary  \ref{corestimate}, we find the lower bound for the eigenvalues of the Dirichlet Laplacian proved by Reilly in \cite{Reilly}. Indeed, on a manifold $M$ with boundary $N$ such that ${\rm Ric}^M \geq k$ with nonnegative mean curvature $H,$  consider any closed 1-form $\alpha$ on $M.$ Take a number $\tau$ big enough so that $\tau\geq \frac{k}{(n-1)H_{min}}$ and $II+\tau\geq 0.$ Then one deduces that $\lambda(\tau,\alpha)\geq \frac{n}{n-1} k.$ As the spectrum of the Robin Laplacian tends to the Dirichlet one when $\tau\to\infty,$ the result then follows.\\

\noindent 
In the last part of this paper, we present two-sided estimates of all the eigenvalues $\lambda_k(\tau,\alpha)$ in terms of $\lambda_1(\tau)=\lambda_1(\tau,0)$ and the Neumann eigenvalues $\lambda_k^N(\alpha):=\lambda_k(0,\alpha),$ using a variational argument (see Theorem \ref{thm:comp} below). These estimates  yield a quantitative  
measurement of the diamagnetism (i.e. the quantity $\lambda(\tau,\alpha)-\lambda_1(\alpha)$). To state this theorem, we define for a normalized eigenfunction of the Robin Laplacian (without magnetic field) $f_\tau:M\to\R$ the constant 
the following constant
\begin{equation}\label{eq:C(tau)}
C(\tau)=\cfrac{\displaystyle\min_{x\in M}f_\tau^2(x) }{\displaystyle\max_{x\in M}f_\tau^2(x) }>0\,.
\end{equation}
Note that $C(0)=1$, $\displaystyle\lim_{\tau\to+\infty}C(\tau)=0$ and the function $f_\tau$ can be selected in a unique manner so that $f_\tau>0$. We have 

\begin{thm}\label{thm:comp}
For all $\tau>0$ and $k\geq 1$, 
$$\lambda_1(\tau)+C(\tau)\lambda_k^N(\alpha)\leq \lambda_k(\tau,\alpha)\leq \lambda_1(\tau)+\frac{1}{C(\tau)}\lambda_k^N(\alpha)\,.$$
\end{thm}

\begin{rem}\label{rem:comp}~
\begin{enumerate}
\item Using the existing estimates on the Neumann eigenvalues $\lambda_k^N(\alpha)$ (see e.g. \cite{CS}), we deduce immediately estimates on the Robin eigenvalues $\lambda_k(\tau,\alpha)$.
\item  {\bf (Zero magnetic field)} Assume that $\alpha$ is closed and not exact. Combining the result in \cite{Sh} and the estimates in Theorem~\ref{thm:comp}, we deduce that $\lambda_1(\tau,\alpha)=\lambda_1(\tau)$ if and only if the flux of $\alpha$ satsifies
$$\Phi_c^\alpha:=\oint_c\alpha \in\mathbb Z$$
for every closed curve $c\subset M$.
\end{enumerate}
\end{rem}

\noindent The rest of the paper is organized as follows. Section~\ref{sec:proofbochner} is devoted to the lengthy proof of Theorem~\ref{bochnermagnetic1}. In Section~\ref{sec:3}, we prove Theorem \ref{estimate1}. Finally, we present the proof of Theorem \ref{thm:comp} in Section~\ref{sec:4}.

\section{Proof of Theorem \ref{bochnermagnetic1}} \label{sec:proofbochner}
In this section, we will prove Theorem \ref{bochnermagnetic1}. We will integrate all the terms in the Bochner formula. First, with the help of the Stokes formula the integral of the l.h.s. of Equation \eqref{eq:bochner} is equal to
$$-\frac{1}{2}\int_M\Delta^M(|d^\alpha f|^2)dv_g=-\frac{1}{2}\int_Ng(d^M(|d^\alpha f|^2),\nu))dv_g=-\int_N\Re\langle \nabla^M_\nu d^\alpha f,d^\alpha f\rangle dv_g.$$
%\begin{eqnarray*}
%-\frac{1}{2}\int_M\Delta^M(|d^\alpha f|^2)dv_g
%&=&\frac{1}{2}\int_M{\rm div}^M(d^M(|d^\alpha f|^2))dv_g\\
%&=&-\frac{1}{2}\int_Ng(d^M(|d^\alpha f|^2),\nu))dv_g\\
%&=&-\int_N\Re\langle \nabla^M_\nu d^\alpha f,d^\alpha f\rangle dv_g.
%\end{eqnarray*}
Now, we will compute the term $\Re\langle \nabla^M_\nu d^\alpha f,d^\alpha f\rangle$ pointwise by decomposing the vectors into the tangential and normal parts over a local orthonormal frame $\{e_i\}_{i=1,\cdots,n-1}$ of $T_xN$ at some point $x\in N.$ Indeed, using the definition of the operator $d^\alpha,$ we write  
\begin{align*}
\langle \nabla^M_\nu d^\alpha f,d^\alpha f\rangle
&=\sum_{i=1}^{n-1}(\nabla^M_\nu d^\alpha f)(e_i) \langle e_i,d^\alpha f\rangle+(\nabla^M_\nu d^\alpha f)(\nu)\langle \nu,d^\alpha f\rangle\\
&=\sum_{i=1}^{n-1}(\nabla^M_\nu d^M f)(e_i) \langle e_i,d^\alpha f\rangle+i\nu(f)\sum_{i=1}^{n-1}\alpha(e_i)\langle e_i,d^\alpha f\rangle\\&\quad+if\sum_{i=1}^{n-1}(\nabla^M_\nu\alpha)(e_i)\langle e_i,d^\alpha f\rangle+(\nabla^M_\nu d^\alpha f)(\nu)\langle \nu,d^\alpha f\rangle\\
&=\sum_{i=1}^{n-1}(\nabla^M_{e_i}d^M f)(\nu) \langle e_i,d^\alpha f\rangle+i\nu(f)\sum_{i=1}^{n-1}\alpha(e_i)\langle e_i,d^\alpha f\rangle\\&\quad+if\sum_{i=1}^{n-1}(d^M\alpha)(\nu,e_i)\langle e_i,d^\alpha f\rangle+if\sum_{i=1}^{n-1}(\nabla^M_{e_i}\alpha)(\nu)\langle e_i,d^\alpha f\rangle\\&\quad+(\nabla^M_\nu d^\alpha f)(\nu)\langle \nu,d^\alpha f\rangle.
\end{align*}
In the last equality, we just use the fact that the hessian of the function $f$ is a symmetric 2-tensor. We then proceed
\begin{align*}
\langle \nabla^M_\nu d^\alpha f,d^\alpha f\rangle
&=\sum_{i=1}^{n-1}e_i(\nu(f))\langle e_i,d^\alpha f\rangle-\sum_{i=1}^{n-1}(d^Mf)(\nabla^M_{e_i}\nu)\langle e_i,d^\alpha f\rangle+i\nu(f)\sum_{i=1}^{n-1}\alpha(e_i)\langle e_i,d^\alpha f\rangle\\
&\quad+if\sum_{i=1}^{n-1}(d^M\alpha)(\nu,e_i)\langle e_i,d^\alpha f\rangle+if\sum_{i=1}^{n-1}e_i(\alpha(\nu))\langle e_i,d^\alpha f\rangle\\
&\quad-if\sum_{i=1}^{n-1}\alpha(\nabla^M_{e_i}\nu)\langle e_i,d^\alpha f\rangle+(\nabla^M_\nu d^\alpha f)(\nu)\langle \nu,d^\alpha f\rangle\\
&=\langle d^N(\nu(f)),d^\alpha f\rangle+\sum_{i=1}^{n-1}(d^Mf)(II(e_i))\langle e_i,d^\alpha f\rangle+i\nu(f)\sum_{i=1}^{n-1}\alpha(e_i)\langle e_i,d^\alpha f\rangle\\&
\quad+if\sum_{i=1}^{n-1}(d^M\alpha)(\nu,e_i)\langle e_i,d^\alpha f\rangle+if\langle d^N(\alpha(\nu)),d^\alpha f\rangle\\
&\quad+if\sum_{i=1}^{n-1}\alpha(II(e_i))\langle e_i,d^\alpha f\rangle+(\nabla^M_\nu d^\alpha f)(\nu)\langle \nu,d^\alpha f\rangle.
\end{align*}
As $\alpha$ is a $1$-form on $M,$ we can write it at any point of the boundary as $\alpha=\alpha^T+\alpha(\nu)\nu.$ We then define the operator $d^\alpha_N$ by $d^\alpha_N h:=d^Nh+ih\alpha^T$ for any complex-valued function $h\in C^\infty(N,\mathbb{C}).$ Hence, the above equality becomes
\begin{align*}
\langle \nabla^M_\nu d^\alpha f,d^\alpha f\rangle=&\langle d^\alpha_N(\nu(f)),d^\alpha f\rangle+\langle II(d^\alpha_Nf),d^\alpha f\rangle+if\langle\nu\lrcorner d^M\alpha,d^\alpha f\rangle\\&+if\langle d^N(\alpha(\nu)),d^\alpha f\rangle+(\nabla^M_\nu d^\alpha f)(\nu)\langle \nu,d^\alpha f\rangle.
\end{align*}
Therefore after integrating, we deduce that 
\begin{align} \label{eq:integralelaplacien}
-\frac{1}{2}\int_M\Delta^M(|d^\alpha f|^2)dv_g=&-\int_N \Re(\langle d^\alpha_N(\nu(f)),d^\alpha f\rangle+\langle II(d^\alpha_Nf),d^\alpha f\rangle+if\langle\nu\lrcorner d^M\alpha,d^\alpha f\rangle\nonumber\\&+if\langle d^N(\alpha(\nu)),d^\alpha f\rangle+(\nabla^M_\nu d^\alpha f)(\nu)\langle \nu,d^\alpha f\rangle)dv_g.
\end{align}
In the second step, we want to integrate the term $\Re\langle d^\alpha f,d^\alpha(\Delta^\alpha f)\rangle$ in the r.h.s. of Theorem \ref{bochnermagnetic}. First, recall the Stokes formula on complex functions: For all $h\in C^\infty(M,\mathbb{C})$ and smooth complex valued $1$-form $\beta,$ one has
$$\int_M \langle d^M h,\beta\rangle dv_g=\int_M h\overline{\delta^M\beta} dv_g-\int_N h\langle\nu,\beta\rangle dv_g.$$
Therefore according to this formula, one can easily get that 
$$\int_M \langle d^\alpha h,\beta\rangle dv_g=\int_M h\overline{\delta^\alpha\beta} dv_g-\int_N h\langle\nu,\beta\rangle dv_g,$$
where the adjoint $\delta^\alpha$ of $d^\alpha$ is given by $\delta^\alpha=\delta^M-i\langle\cdot,\alpha\rangle$ \cite[Def. 2.1]{ELMP}. Here we mention that $\delta^\alpha X=-{\rm trace}(\nabla^\alpha X),$ where $\nabla^\alpha$ is the magnetic covariant derivative defined previously. Hence, by taking $h=\Delta^\alpha f$ and $\beta=d^\alpha f,$ we deduce
\begin{eqnarray}\label{eq:integraledalpha}
\int_M\langle d^\alpha(\Delta^\alpha f),d^\alpha f\rangle dv_g=\int_M |\Delta^\alpha f|^2dv_g-\int_N(\Delta^\alpha f)\langle\nu,d^\alpha f\rangle dv_g.
\end{eqnarray} 
Now we want to evaluate the term $\Delta^\alpha f$ in the second integral of the r.h.s. of the equality above. Using the compatibility equations in \cite[Lem. 3.2]{ELMP} and taking an orthonormal frame $\{e_i\}_{i=1,\cdots,n-1}$ of $TN$ with $\nabla^N_{e_i}e_i=0$ at some point, we compute
\begin{align*} 
\Delta^\alpha f&=-\sum_{i=1}^{n-1}\langle\nabla^\alpha_{e_i}(d^\alpha f),e_i\rangle-\langle\nabla^\alpha_{\nu}(d^\alpha f),\nu\rangle\\
&=-\sum_{i=1}^{n-1}e_i(\langle d^\alpha f,e_i\rangle)+\sum_{i=1}^{n-1}\langle d^\alpha f,\nabla^\alpha_{e_i}e_i\rangle-\langle\nabla^\alpha_{\nu}(d^\alpha f),\nu\rangle\\
&=-\sum_{i=1}^{n-1}e_i(\langle d^\alpha f,e_i\rangle)+\sum_{i=1}^{n-1}\langle d^\alpha f,\nabla^M_{e_i}e_i+i\alpha(e_i)e_i\rangle-\langle \nabla^\alpha_{\nu}(d^\alpha f),\nu\rangle\\
&=-\sum_{i=1}^{n-1}e_i(\langle d^\alpha f,e_i\rangle)+\sum_{i=1}^{n-1}\langle d^\alpha f,II(e_i,e_i)\nu+i\alpha(e_i)e_i\rangle-\langle\nabla^\alpha_{\nu}(d^\alpha f),\nu\rangle\\
&=-\sum_{i=1}^{n-1}e_i(\langle d_N^\alpha f,e_i\rangle)+(n-1)H \langle d^\alpha f,\nu\rangle+\sum_{i=1}^{n-1}\langle d_N^\alpha f,i\alpha(e_i)e_i\rangle-\langle\nabla^\alpha_{\nu}(d^\alpha f),\nu\rangle\\
&=\Delta_N^\alpha f+(n-1)H \langle d^\alpha f,\nu\rangle-\langle\nabla^\alpha_{\nu}(d^\alpha f),\nu\rangle, 
\end{align*}
where $\Delta_N^\alpha :=\delta_N^\alpha d_N^\alpha,$ with $\delta_N^\alpha=\delta^N-i(\cdot,\alpha^T).$ We notice that $\delta_N^\alpha$ is the $L^2$-adjoint of $d_N^\alpha$ on $N.$ Plugging the expression of $\Delta^\alpha f$ above into Equation $\eqref{eq:integraledalpha},$ we find 
\begin{eqnarray}\label{dalphadelta}
\int_M\langle d^\alpha(\Delta^\alpha f),d^\alpha f\rangle dv_g=\int_M |\Delta^\alpha f|^2dv_g-\int_N(\Delta_N^\alpha f)\langle\nu,d^\alpha f\rangle dv_g-(n-1)\int_N H|\langle d^\alpha f,\nu\rangle|^2dv_g\nonumber\\+\int_N\langle\nabla^\alpha_{\nu}(d^\alpha f),\nu\rangle\langle\nu,d^\alpha f\rangle dv_g.\nonumber\\
=\int_M |\Delta^\alpha f|^2dv_g-\int_N(\Delta_N^\alpha f)\langle\nu,d^\alpha f\rangle dv_g-(n-1)\int_N H|\langle d^\alpha f,\nu\rangle|^2dv_g\nonumber\\+\int_N\langle\nabla^M_{\nu}(d^\alpha f),\nu\rangle\langle\nu,d^\alpha f\rangle dv_g+\int_N i\alpha(\nu)|\langle\nu,d^\alpha f\rangle|^2dv_g.\nonumber\\
\end{eqnarray} 
The last step is to compute the term $\frac{i}{2}\displaystyle\int_M\langle \bar{f}d^\alpha f,\delta^M d^M\alpha\rangle dv_g$ and its conjugate in Theorem \ref{bochnermagnetic}. For this, we proceed as in \cite[p.17]{ELMP} to get
\begin{eqnarray} \label{eq:barfdalpha}
\frac{i}{2}\int_M\langle \bar{f}d^\alpha f,\delta^M d^M\alpha\rangle dv_g&=&\frac{i}{2}\int_M \langle d^M(\bar{f}d^\alpha f), d^M\alpha\rangle dv_g+\frac{i}{2}\int_N\langle \bar{f}d^\alpha f, \nu\lrcorner d^M\alpha\rangle dv_g\nonumber\\
&=&\frac{i}{2}\int_M (d^M\alpha)(\overline{d^\alpha f},d^\alpha f) dv_g-\frac{1}{2}\int_M |f|^2|d^M\alpha|^2 dv_g\nonumber\\&&+\frac{i}{2}\int_N\langle \bar{f}d^\alpha f, \nu\lrcorner d^M\alpha\rangle dv_g\nonumber.\\
\end{eqnarray} 
Now, we have all the ingredients to integrate Equation \eqref{eq:bochner} over $M.$ In fact, using Equations \eqref{eq:integralelaplacien}, \eqref{dalphadelta} and \eqref{eq:barfdalpha}, we find that
\begin{eqnarray*} 
-\int_N\Re(\langle d^\alpha_N(\nu(f)),d^\alpha f\rangle+\langle II(d^\alpha_Nf),d^\alpha f\rangle+if\langle\nu\lrcorner d^M\alpha,d^\alpha f\rangle
+if\langle d^N(\alpha(\nu)),d^\alpha f\rangle\\
+(\nabla^M_\nu d^\alpha f)(\nu)\langle \nu,d^\alpha f\rangle)dv_g=\int_M|{\rm Hess}^\alpha f|^2dv_g-\int_M |\Delta^\alpha f|^2dv_g+\int_N\Re((\Delta_N^\alpha f)\langle\nu,d^\alpha f\rangle) dv_g\\
+(n-1)\int_N H|\langle d^\alpha f,\nu\rangle|^2dv_g-\int_N\Re(\langle\nabla^M_{\nu}(d^\alpha f),\nu\rangle\langle\nu,d^\alpha f\rangle) dv_g+\int_M {\rm Ric}^M(d^\alpha f,d^\alpha f) dv_g\\
+\frac{i}{2}\int_M \left(\underbrace{(d^M\alpha)(d^\alpha f,\overline{d^\alpha f})-(d^M\alpha)(\overline{d^\alpha f},d^\alpha f)}_{2i\Im m((d^M\alpha)(d^\alpha f,\overline{d^\alpha f}))}\right) dv_g
-\int_M|f|^2|d^M\alpha|^2dv_g\\
+\frac{i}{2}\int_N \left(\underbrace{\langle \bar{f}d^\alpha f, \nu\lrcorner d^M\alpha\rangle-\langle f\overline{d^\alpha f}, \nu\lrcorner d^M\alpha}_{=-2i \Im m f\langle \nu\lrcorner d^M\alpha,d^\alpha f \rangle}\rangle\right) dv_g.
\end{eqnarray*}

\noindent By writing $d^\alpha f=d_N^\alpha f+(\nu(f)+if\alpha(\nu))\nu$ at any point of the boundary, the first integral in the l.h.s. reduces to 
\begin{eqnarray*} 
\int_N\Re\langle d^\alpha_N(\nu(f)),d^\alpha f\rangle dv_g&=&\int_N\Re\langle d^\alpha_N(\nu(f)),d_N^\alpha f\rangle dv_g\\
&=&\int_N\Re (\nu(f)\overline{\delta_N^\alpha d_N^\alpha f}) dv_g=\int_N\Re (\nu(f)\overline{\Delta_N^\alpha f}) dv_g\\
&=&\int_N\Re (\langle d^\alpha f-i\alpha f,\nu\rangle\overline{\Delta_N^\alpha f}) dv_g\\
%&=&\int_N\Re (\langle d^\alpha f,\nu\rangle\overline{\Delta_N^\alpha f})dv_g-\int_N\Re(i \alpha(\nu)f\, \overline{\Delta_N^\alpha f}) dv_g\\
&=&\int_N\Re (\langle\nu,d^\alpha f\rangle\Delta_N^\alpha f)dv_g-\int_N\Re(i \alpha(\nu)f\, \overline{\Delta_N^\alpha f}) dv_g.
\end{eqnarray*}
Using the fact that $\delta_N^\alpha$ is the $L^2$-adjoint of $d_N^\alpha$ and that $d_N^\alpha(f_1f_2)=f_2d^Nf_1+f_1d_N^\alpha f_2$ for any complex valued functions $f_1$ and $f_2$ on $N,$ the above equality becomes
\begin{eqnarray*} 
\int_N\Re\langle d^\alpha_N(\nu(f)),d^\alpha f\rangle dv_g&=&\int_N\Re (\langle\nu,d^\alpha f\rangle\Delta_N^\alpha f)dv_g-\int_N\Re\langle d_N^\alpha f,d_N^\alpha\left(i\alpha(\nu) f\right)\rangle dv_g\\
&=&\int_N\Re (\langle\nu,d^\alpha f\rangle\Delta_N^\alpha f)dv_g+\int_N\Re (i\langle d^\alpha_Nf,fd^N(\alpha(\nu))+\alpha(\nu)d_N^\alpha f\rangle) dv_g\\
&=&\int_N\Re (\langle\nu,d^\alpha f\rangle\Delta_N^\alpha f)dv_g+\int_N\Re (i\bar{f}\langle d^\alpha_Nf,d^N(\alpha(\nu))\rangle) dv_g\\&&+\int_N\alpha(\nu)\underbrace{\Re (i\langle d^\alpha_Nf,d_N^\alpha f\rangle)}_{=0} dv_g\\
&=&\int_N\Re (\langle\nu,d^\alpha f\rangle\Delta_N^\alpha f)dv_g-\int_N\Re (i f\langle d^N(\alpha(\nu)),d^\alpha f\rangle) dv_g.
\end{eqnarray*}
Therefore, we deduce\\\\
%\begin{eqnarray*}
$$-2\int_N\Re (\langle\nu,d^\alpha f\rangle\Delta_N^\alpha f)dv_g-\int_N\langle II(d^\alpha_Nf),d_N^\alpha f\rangle dv_g=$$
$$\int_M|{\rm Hess}^\alpha f|^2dv_g-\int_M |\Delta^\alpha f|^2dv_g+(n-1)\int_N H|\langle d^\alpha f,\nu\rangle|^2dv_g+\int_M {\rm Ric}^M(d^\alpha f,d^\alpha f) dv_g$$
$$-\int_M\Im m\left((d^M\alpha)(d^\alpha f,\overline{d^\alpha f})\right)dv_g
-\int_M|f|^2|d^M\alpha|^2dv_g.$$
The proof of the proposition then follows.
\hfill$\square$

\section{Proof of Theorem \ref{estimate1}}\label{sec:3}
In the following, we will give a proof of Theorem \ref{estimate1}. For this, we consider an eigenfunction $f$ of the Robin Laplacian associated to the eigenvalue $\lambda(\tau,\alpha),$ that is $\Delta^\alpha f=\lambda(\tau,\alpha) f$ with $\nu(f)+if\alpha(\nu)=\tau f$ for some positive $\tau.$ We then apply Equality \eqref{bochnermagnetic2} to the eigenfunction $f$. First, we have
$$\int_N\Re (\langle\nu,d^\alpha f\rangle\Delta_N^\alpha f)dv_g= \tau\int_N\Re (\bar{f}\Delta_N^\alpha f)dv_g=\tau\int_N\Re (f\overline{\Delta_N^\alpha f})dv_g=\tau\int_N|d^\alpha_N f|^2dv_g.$$
Also, the following inequality 
$$\int_M \Im m\left((d^M\alpha)(d^\alpha f,\overline{d^\alpha f})\right)dv_g\leq ||d^M\alpha||_\infty\int_M|d^\alpha f|^2dv_g,$$ 
holds. Therefore, as the r.h.s. of Equality \eqref{bochnermagnetic2} is nonnegative, we get after using the conditions ${\rm Ric}^M\geq k$ and $II+\tau\geq 0$ that 
\begin{eqnarray*}
0&\leq &\frac{n-1}{n}\lambda(\tau,\alpha)^2\int_M|f|^2dv_g-(k-||d^M\alpha||_\infty)\int_M|d^\alpha f|^2dv_g+||d^M\alpha||_\infty^2\int_M|f|^2dv_g\\&&-(n-1)\tau^2\int_N H|f|^2dv_g-\tau\int_N|d^\alpha_N f|^2dv_g.
\end{eqnarray*}
Since $f$ is an eigenfunction of the Laplacian, one has
\begin{equation*}
\int_M |d^\alpha f|^2dv_g=\lambda(\tau,\alpha)\int_M |f|^2dv_g-\tau\int_N|f|^2 dv_g.
\end{equation*} 
Hence, the above inequality reduces to 
\begin{eqnarray*}
0&\leq &\frac{n-1}{n}\lambda(\tau,\alpha)^2\int_M|f|^2dv_g-(k-||d^M\alpha||_\infty)\lambda(\tau,\alpha)\int_M| f|^2dv_g+(k-||d^M\alpha||_\infty)\tau\int_N|f|^2dv_g\\&&+||d^M\alpha||_\infty^2\int_M|f|^2dv_g-(n-1)\tau^2H_{\min}\int_N |f|^2dv_g-\tau\int_N|d^\alpha_N f|^2dv_g.
\end{eqnarray*}
By grouping the terms and using the fact that the last term is nonpositive, we find at the end 
\begin{eqnarray*}
0&\leq &\left(\frac{n-1}{n}\lambda(\tau,\alpha)^2-(k-||d^M\alpha||_\infty)\lambda(\tau,\alpha)+||d^M\alpha||_\infty^2\right)\int_M|f|^2dv_g\\&&+\tau\left(k-||d^M\alpha||_\infty-(n-1)\tau H_{\min}\right)\int_N|f|^2dv_g.
\end{eqnarray*}
Since now the sign of the term $(k-||d^M\alpha||_\infty)-(n-1)\tau H_{\min}$ is nonpositive, we deduce as in \cite[Eq. 62]{ELMP} the inequality
$$0\leq \frac{n-1}{n}\lambda(\tau,\alpha)^2-(k-||d^M\alpha||_\infty)\lambda(\tau,\alpha)+||d^M\alpha||_\infty^2.$$ 
Therefore, as the discriminant of this polynomial is nonnegative, we finish the proof. 
\hfill$\square$

\section{Proof of Theorem~\ref{thm:comp}}\label{sec:4}
Let $f$ be the function defined by $f=u f_\tau$, where $u:M\rightarrow \C$ is a complex valued function on $M$ and $f_\tau$ is a normalized eigenfunction of the Robin Laplacian associated to the first eigenvalue $\lambda_1(\tau).$ Then, we compute
\begin{eqnarray*}
\int_M|(d^M+i\alpha)f|^2dv_g&=&\int_\Omega|ud^Mf_\tau+f_\tau(d^Mu+i\alpha u)|^2dv_g\\
&=&\int_M |u|^2|d^Mf_\tau|^2dv_g+\int_M f_\tau^2|(d^M+i\alpha)u|^2dv_g\\&&+2\int_M f_\tau\Re\langle ud^Mf_\tau,d^Mu+i\alpha u\rangle dv_g\\
&=&\int_M f_\tau\delta^M(|u|^2d^Mf_\tau)dv_g-\tau\int_{N}|u|^2f_\tau^2dv_g+\int_M f_\tau^2|(d^M+i\alpha)u|^2dv_g\\&&+\int_M \Re\langle d^M(f_\tau^2),{\bar u}d^Mu\rangle dv_g\\
&=&\int_M f_\tau|u|^2\delta^M(d^Mf_\tau)dv_g-\int_M f_\tau g(d^M(|u|^2),d^M(f_\tau)) dv_g-\tau\int_{N}|u|^2f_\tau^2dv_g\\&&+\int_M f_\tau^2|(d^M+i\alpha)u|^2+\int_M \Re\langle d^M(f_\tau^2),{\bar u}d^Mu\rangle dv_g\\
&=&\lambda_1(\tau) \int_M f_\tau^2|u|^2dv_g-\int_M f_\tau g(d^M(|u|^2),d^M(f_\tau)) dv_g-\tau\int_{N}|u|^2f_\tau^2dv_g\\&&+\int_M f_\tau^2|(d^M+i\alpha)u|^2dv_g+\int_M \Re\langle d^M(f_\tau^2),{\bar u}d^Mu\rangle dv_g.
\end{eqnarray*}
Now, it is easy to see that one has pointwise
$$f_\tau g(d^M(|u|^2),d^M(f_\tau))=f_\tau\langle {\bar u}d^Mu+ud^M\overline{u},d^M(f_\tau)\rangle=\Re\langle d^M(f_\tau^2),{\bar u}d^Mu\rangle.$$
Consequently, we deduce that 
$$\frac{\int_M|d^\alpha f |^2dv_g+\tau\int_{N}f ^2dv_g}{||f||^2}=\lambda_1(\tau)+\frac{\int_M f_\tau^2 |d^\alpha u|^2\,dv_g}{\int_M |u|^2f_\tau^2\,dv_g} \,.$$
Now the proof follows from the variational min-max principle.
Indeed,  the definition of $C(\tau)$ in \eqref{eq:C(tau)} yields
$$C(\tau)\frac{\int_M  |d^\alpha u|^2\,dv_g}{\int_M |u|^2\,dv_g} \leq \frac{\int_M f_\tau^2 |d^\alpha u|^2\,dv_g}{\int_M |u|^2f_\tau^2\,dv_g} \leq \frac1{C(\tau)} \frac{\int_M  |d^\alpha u|^2\,dv_g}{\int_M |u|^2\,dv_g},$$
which finishes the proof.
\hfill$\square$\\

\noindent{\bf Acknowledgment.} The authors are indebted to Nicolas Ginoux, Norbert Peyerimhoff and Alessandro Savo for their valuable comments on the paper. The first named author acknowledges the financial support of the Alexander von Humboldt Foundation. The authors are indebted to the anonymous referee who suggested the proof of Corollary~\ref{corestimate}.

\end{document}